\RequirePackage{fix-cm}
\documentclass[smallcondensed]{svjour3}     
\smartqed  
\usepackage[parfill]{parskip}    
\usepackage{amsmath, amssymb, latexsym, enumerate, mathrsfs}
\usepackage{graphicx}
\usepackage{epstopdf}
\usepackage{subfig}

\usepackage{hyperref}
\usepackage{array}
\usepackage{times}
\usepackage{moreverb}
\usepackage{marvosym}
\journalname{}
\begin{document}

\title{On the risk-sensitive escape control for diffusion processes pertaining to an expanding construction of distributed control systems
}

\titlerunning{On the risk-sensitive escape control for diffusion processes}        

\author{Getachew K. Befekadu \and Panos~J.~Antsaklis
}

\institute{G. K. Befekadu~ ({\large\Letter}\negthinspace) \at
	          Department of Electrical Engineering, University of Notre Dame, Notre Dame, IN 46556, USA. \\
	          Tel.: +1 574 631 6618\\
	          Fax: +1 574 631 4393 \\
	          \email{gbefekadu1@nd.edu}           
	          \and
	          P. J. Antsaklis \at
	          Department of Electrical Engineering, University of Notre Dame, Notre Dame, IN 46556, USA. \\
	          \email{antsaklis.1@nd.edu}           
}

\date{Received: September 21, 2014 / Accepted: date}

\maketitle

\begin{abstract}
In this paper, we consider an expanding construction of a distributed control system, which is obtained by adding a new subsystem one after the other, until all $n$ subsystems, where $n \ge 2$, are included in the distributed control system. It is assumed that a small random perturbation enters only into the first subsystem and is then subsequently transmitted to the other subsystems. Moreover, for any $\ell \in \{2, \ldots, n\}$, the distributed control system, compatible with the expanding construction, which is obtained from the first $\ell$ subsystems, satisfies an appropriate H\"{o}rmander condition. As a result of this, the diffusion process is degenerate, i.e., the backward operator associated with it is a degenerate parabolic equation. Our main interest here is to prevent the diffusion process (that corresponds to a particular subsystem) from leaving a given bounded open domain. In particular, we consider a risk-sensitive version of the mean escape time criterion with respect to each of the subsystems. Using a variational representation, we characterize the risk-sensitive escape control for the diffusion process as the lower and upper values of an associated stochastic differential game. Finally, we comment on the implication of our results, where one is also interested in evaluating the performance of the risk-sensitive escape control, when there is some modeling error in the distributed control system.

\end{abstract}

\section{Introduction}	 \label{S1}
\begin{figure}[tbh]
\vspace{-3mm}
\begin{center}
\hspace{-10mm}\subfloat{\includegraphics[width=28mm]{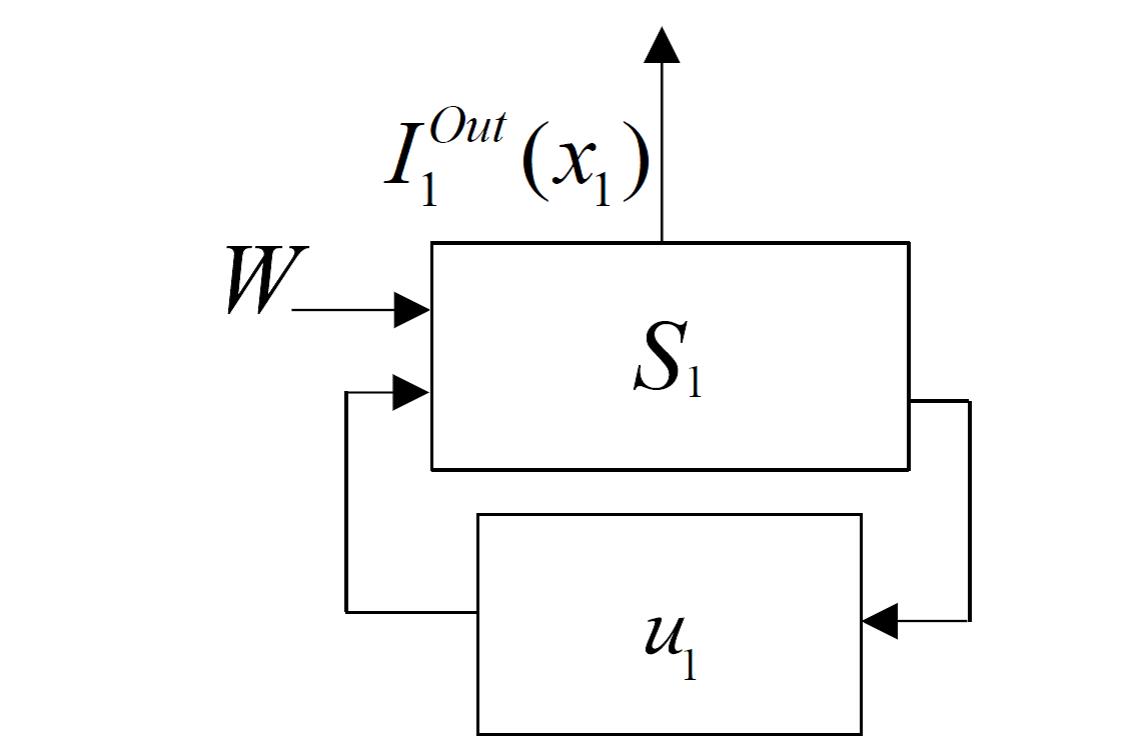}}
\subfloat{\includegraphics[width=24mm]{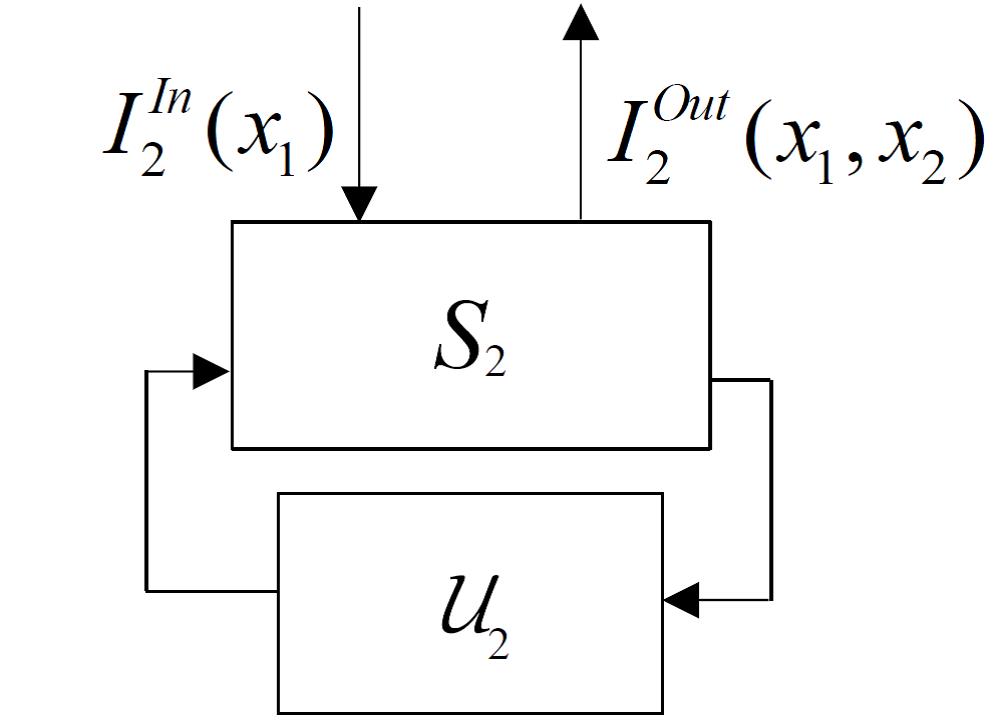}}$~~~~~~~$
\subfloat{\includegraphics[width=34mm]{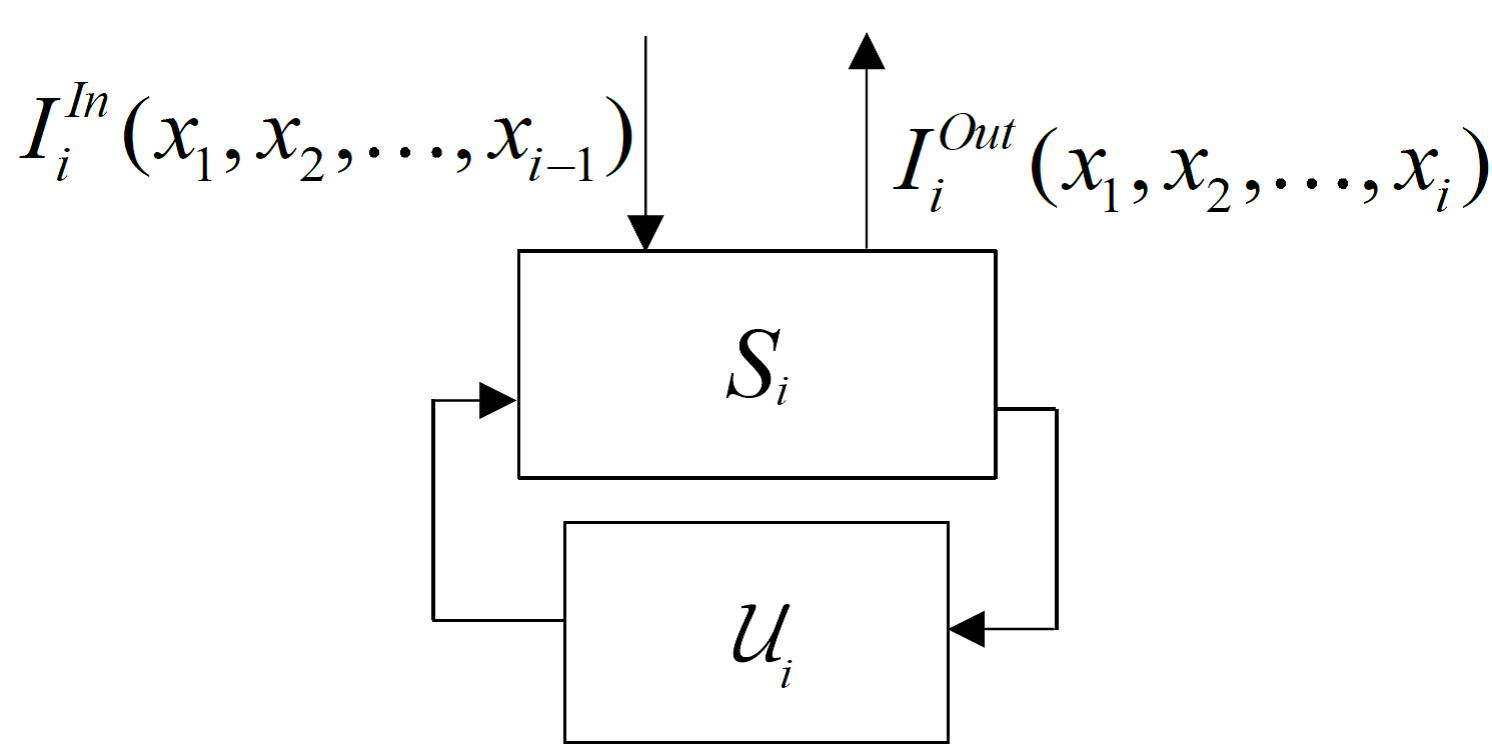}} $~~~~~~~~$
\subfloat{\includegraphics[width=25mm]{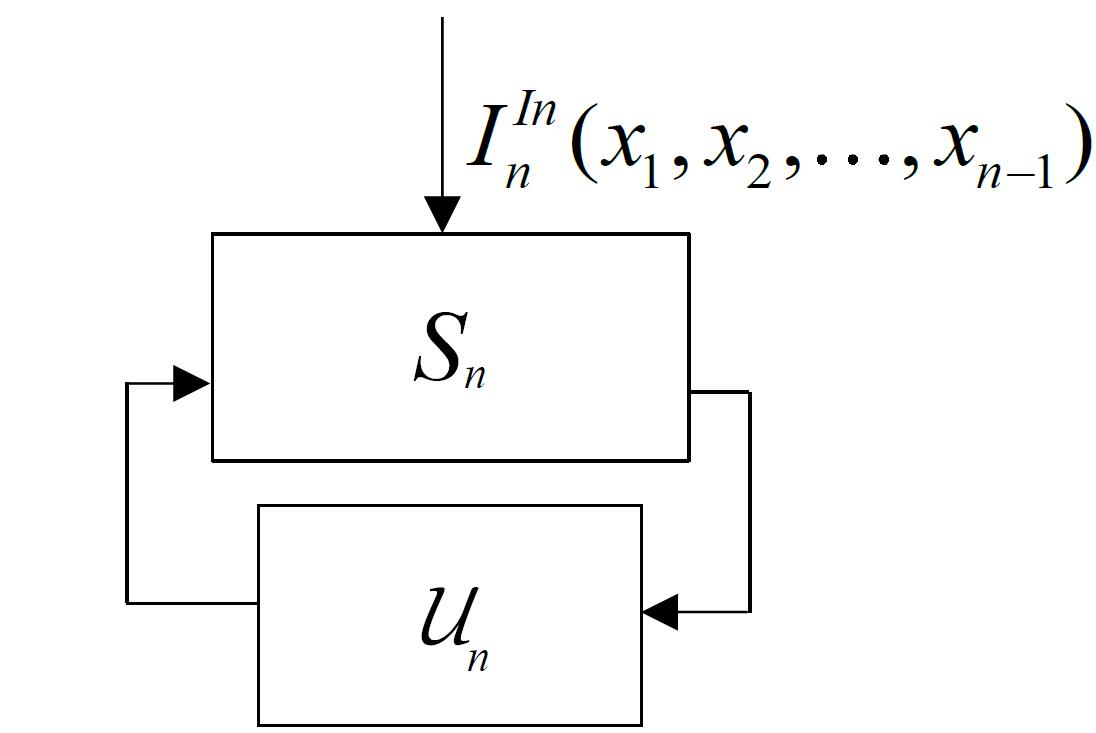}}\\
{$ \begin{array}{l@{\ }c@{\ }l}\\
\text{where} & \\
&S_1:  \, \, d x^{1}(t) = m_1\bigl(x^{1}(t), u_1(t)\bigr) dt + \sqrt{\epsilon} \sigma\bigl(x^{1}(t)\bigr)dW(t), \\
&S_i: \,\, d x^{i}(t) = m_{i}\bigl(x^{1}(t), x^{2}(t), \ldots, x^{i}(t), u_{i}(t)\bigr) dt, ~ i = 2, 3, \ldots n,\\
   &   u_j(t) = \gamma_j^s \bigl(t, W(r), 0 \le r \le t \bigr), ~ t \in [0,\, s], ~ \forall s \ge 0, ~ j =1, 2,  \ldots n,\\
     & I_i^{In}(x^1, x^2, \ldots, x^{i-1}) ~ \text{and} ~  I_i^{Out}(x^1, x^2, \ldots, x^{i}) ~ \text{are information for the expanding construction.}
\end{array}$}
\caption{Distributed control systems with small random perturbations} \label{Fig-DCS}
\vspace{-3 mm}
\end{center}
\end{figure}

We consider the diffusion processes $\bigl(x^{1}(t), x^{2}(t), \ldots, x^{n}(t)\bigr)$ pertaining to the following distributed control system, with small random perturbations (see Fig.~\ref{Fig-DCS})\footnote{This work is, in some sense, a continuation of our previous paper \cite{BefAn14}.}
\begin{align}
\left.\begin{array}{l}
d x^{1}(t) = m_1\bigl(x^{1}(t), u_1(t)\bigr) dt + \sqrt{\epsilon} \sigma\bigl(x^{1}(t)\bigr)dW(t) \\
d x^{2}(t) = m_2\bigl(x^{1}(t), x^{2}(t), u_2(t)\bigr) dt  \\
 \quad\quad\quad~ \vdots  \\
d x^{n}(t) = m_n\bigl(x^{1}(t), x^{2}(t), \ldots, x^{n}(t), u_{n}(t)\bigr) dt \\
 x^{1}(0)=x_0^{1}, \,\, x^{2}(0)=x_0^{2}, \,\, \ldots, \,\, x^{n}(0)=x_0^{n}, \,\, t \ge 0, \,\, n \ge 2
\end{array}\right\}  \label{Eq1} 
\end{align}
where
\begin{itemize}
\item[-] $x^{i}(\cdot)$ is an $\mathbb{R}^{d}$-valued diffusion process that corresponds to the $i$th-subsystem (with $i \in \{1,2, \ldots, n\}$),
\item[-] the functions $m_i \colon \underbrace{\mathbb{R}^d \times \mathbb{R}^d \times \cdots \times \mathbb{R}^d}_{i - \rm{times}} \times \,\mathcal{U}_i \rightarrow \mathbb{R}^{d}$ are uniformly Lipschitz, with bounded first derivatives, $\epsilon$ is a small positive number (which is related to the random perturbation level in the system),
\item[-] $\sigma \colon \mathbb{R}^{d} \rightarrow \mathbb{R}^{d \times m}$ is Lipschitz with the least eigenvalue of $\sigma(\cdot)\,\sigma^T(\cdot)$ uniformly bounded away from zero, i.e., 
\begin{align*}
 \sigma(x)\,\sigma^T(x) \ge \lambda I_{d \times d} , \quad \forall x \in \mathbb{R}^{d},
\end{align*}
for some $\lambda > 0$,
\item[-] $W(\cdot)$ (with $W(0) = 0$) is a $m$-dimensional standard Wiener process,
\item[-] $u_i(\cdot)$ is a $\,\mathcal{U}_i$-valued measurable control process to the $i$th-subsystem, i.e., an admissible control from the measurable set \,$\mathcal{U}_i \subset \mathbb{R}^{r_i}$.
\end{itemize}

In this paper, we identify two admissible controls $u_{i},\, \tilde{u}_i \in \mathcal{U}_{i}$, for $i = 1, 2 \ldots, n$, being the same on $[0,\, s]$ if $\operatorname{Prob}\bigl\{u_i(t) = \tilde{u}_i(t),\, \forall t \in [0,\, s] \bigr\} =1$. If $u_i \in \mathcal{U}_i$, then, for every $s \in [0,\, \infty)$, there exists a Borel measurable function $\gamma_i^s \colon \mathscr{C} \bigl([0,\,T], \mathbb{R}^m \bigr) \rightarrow \mathcal{U}_i$ (with respect to some underlying Borel $\sigma$-algebra) such that 
\begin{align}
u_i(t) = \gamma_i^s \bigl(t, W(r), 0 \le r \le t \bigr), \quad t \in [0,\, s], \label{Eq2}
\end{align}
with probability one (w.p.1).

The functions $m_{\ell}$, for $\ell \in \{2, \ldots, n\}$, in Equation~\eqref{Eq1}, with any progressively measurable control $u_{i}$, depend only on $\bigl(x^{1}, x^{2}, \ldots, x^{\ell}\bigr)$. Furthermore, we assume that the distributed control system, which is formed by the first $\ell$ subsystems, satisfies an appropriate H\"{o}rmander condition, i.e., a hypoellipticity assumption on the diffusion processes $\bigl(x^{1}(t), \\ x^{2}(t), \cdots, x^{\ell}(t)\bigr)$ (e.g., see \cite{Hor67} or \cite[Section~3]{Ell73}). Notice that the random perturbation has to pass through the second subsystem, the third subsystem, and so on to the $\ell$th-subsystem. Hence, such a distributed control system is described by an $\ell \times d$ dimensional diffusion process, which is degenerate in the sense that the backward operator associated with it is a degenerate parabolic equation.

\begin{figure}
\begin{center}
\setlength{\unitlength}{2.55mm}
\begin{picture}(50,13)
\put(1,9.5){\line(1,0){45}}
\put(1,13){\line(0,-1){3.5}}
\put(46,13){\line(0,-1){3.5}}
\put(18,11){{$ 
\begin{array}{r@{\ }c@{\ }l}
{\rm Interconnection}
\end{array}$}}
\put(1,13){\line(1,0){45}}

\put(18.75,4){\line(1,0){0.2}}
\put(19.50,4){\line(1,0){0.2}}
\put(20.25,4){\line(1,0){0.2}}

\put(33.75,4){\line(1,0){0.2}}
\put(34.50,4){\line(1,0){0.2}}
\put(35.25,4){\line(1,0){0.2}}
\end{picture}
\end{center}
\vspace{-28.35mm}
\end{figure}

\begin{remark} \label{R1}
In general, the hypoellipticity is related to a strong accessibility property of controllable nonlinear systems that are driven by white noise (e.g., see \cite{SusJu72} concerning the controllability of nonlinear systems, which is closely related to \cite{StrVa72} and \cite{IchKu74}). That is, the hypoellipticity assumption implies that the diffusion process $x^{\ell}(t)$ has a transition probability density $p_{(t)}^{\ell}\bigl((x^{1}, \ldots, x^{\ell}), \otimes_{i=1}^{\ell} d\mu_i\bigr)$, which is $C^{\infty}$ on $\mathbb{R}^{2(d \times \ell)}$, with a strong Feller property.
\end{remark}
Let $D_i \subset \mathbb{R}^{d}$, for $i =1, 2, \ldots, n$, be bounded open domains with smooth boundaries (i.e., $\partial D_i$ is a manifold of class $C^2$). Moreover, let $\Omega_{\ell}$ be the open sets that are given by
\begin{align*}
\Omega_{\ell} = D_1 \times D_2 \times \cdots \times D_{\ell}, \quad \ell \in \{2, 3, \ldots, n\}.
\end{align*}
Suppose that, for a fixed $\ell \in \{2, 3, \ldots, n\}$, the distributed control system, which is compatible with expanding construction, is formed by the first $\ell$ subsystems (i.e., obtained by adding one after the other, until all $\ell$th subsystems are included). Furthermore, assume that the newly constructed distributed control system is composed with some admissible controls $u_{i}^{\ast}(t) \in \mathcal{U}_i$, $\forall t \in [0,\, \infty)$, for $i =1, 2, \ldots, \ell-1$. Let $\tau_{\ell}^{\epsilon}=\tau_{\ell}^{\epsilon}(x^{1}, x^{2}, \ldots, x^{\ell}, u_{\ell})$ be the exit-time for the diffusion process $x^{\ell}(t)$ (corresponding to the $\ell$th-subsystem), for a fixed $\epsilon > 0$, with $u_{\ell} \in \mathcal{U}_{\ell}$, from the given domain $D_{\ell}$, i.e., 
\begin{align}
\tau_{\ell}^{\epsilon} = \inf \Bigl\{ t > 0 \, \bigl\vert \, x^{\ell}(t) \notin D_{\ell} \Bigr\}, \label{Eq3}
\end{align}
which depends on the behavior of the following (deterministic) distributed control system
\begin{align}
\left.\begin{array}{l}
d \xi^{1}(t) = m_1\bigl(\xi^{1}(t), u_1^{\ast}(t)\bigr) dt \\
d \xi^{2}(t) = m_2\bigl(\xi^{1}(t), \xi^{2}(t), u_2^{\ast}(t)\bigr) dt  \\
 \quad\quad\quad~ \vdots  \\
d \xi^{\ell-1}(t) = m_{\ell-1}\bigl(\xi^{1}(t), \xi^{2}(t), \ldots, \xi^{\ell-1}(t), u_{\ell-1}^{\ast}(t)\bigr) dt\\ 
d \xi^{\ell}(t) = m_{\ell}\bigl(\xi^{1}(t), \xi^{2}(t), \ldots, \xi^{\ell}(t), u_{\ell}(t)\bigr) dt\\
  \bigl(\xi^{1}(0), \xi^{2}(0), \ldots, \xi^{\ell}(0)\bigr) \triangleq \bigl(x_0^{1}, x_0^{2}, \ldots, x_0^{\ell}\bigr) \in \Omega_{\ell}, \,\, t \ge 0
\end{array}\right\}   \label{Eq4}
\end{align}
In this paper, we specifically consider a risk-sensitive version of the mean escape time criterion with respect to the $\ell$th-subsystem, i.e.,
\begin{align}
-\epsilon \log \mathbb{E}_{x_{\widehat{1,\ell}}}^{\epsilon} \exp \biggl\{-\frac{1}{\epsilon} \theta_{\ell}\, \tau_{\ell}^{\epsilon} \biggr\},  \label{Eq5}
\end{align}
where $\theta_{\ell}$, for each $\ell \in \{2, 3, \ldots, n\}$, are positive design parameters and the expectation $\mathbb{E}_{x_{\widehat{1,\ell}}}^{\epsilon} \bigl\{\cdot\bigr\}$ is conditioned on the initial point $\bigl(x_0^{1}, x_0^{2}, \ldots, x_0^{\ell}\bigr) \in \Omega_{\ell}$ as well as on the admissible controls $\bigl(u_{1}^{\ast}, u_2^{\ast}, \ldots, u_{\ell-1}^{\ast}, u_{\ell}\bigr) \in \prod\nolimits_{i=1}^{\ell} \mathcal{U}_i$. Notice that $\tau_{1}^{\epsilon}$ in the exit-time for the diffusion process $x^1(t)$ (which corresponds to the $1$st-subsystem) from the domain $D_{1}$ with respect to the admissible (optimal) control $u_{1}(t) \in \mathcal{U}_1$, $\forall t \in [0,\,\infty)$, with $\theta_{1} > 0$.\footnote{$\mathbb{E}_{x_{\widehat{1,\ell}}}^{\epsilon} \bigl\{\cdot \bigr\} \triangleq \mathbb{E}_{x^{1}, x^{2}, \ldots, x^{\ell}} \bigl\{\cdot\bigr\}$}

\begin{remark} \label{R2}
Here we remark that the criterion in Equation~\eqref{Eq5} makes sense only if we have the following conditions
\begin{align}
\tau_{1}^{\epsilon}\, \ge\, \tau_{2}^{\epsilon}\, \ge \,\cdots \,\ge \, \tau_{\ell}^{\epsilon}.  \label{Eq6}
\end{align}
Moreover, such conditions depend on the constituting subsystems, the admissible controls from the measurable sets $\prod_{i=1}^{\ell} \mathcal{U}_i$, as well as on the given bounded open domains $D_{i}$, for $i = 1, 2, \ldots, \ell$ (see Section~\ref{S3(2)} for further discussion).
\end{remark}

Then, the problem of risk-sensitive escape control (with respect to the $\ell$th-subsystem) will amount to obtaining a supremum value for $V_{\ell, \theta_{\ell}}^{\epsilon} = V_{\ell, \theta_{\ell}}^{\epsilon} \bigl(x_0^{1}, x_0^{2}, \ldots, x_0^{\ell}\bigr)$, i.e.,
\begin{align}
V_{\ell, \theta_{\ell}}^{\epsilon} \triangleq \sup_{u_{\ell} \in \mathcal{U}_{\ell}} -\epsilon \log \mathbb{E}_{x_{\widehat{1,\ell}}}^{\epsilon} \exp \biggl\{-\frac{1}{\epsilon} \theta_{\ell}\, \tau_{\ell}^{\epsilon} \biggr\},   \label{Eq7}
\end{align}
with respect to some progressively measurable control $u_{\ell} \in \mathcal{U}_{\ell}$, for each $\ell \in \{2, 3, \ldots, n\}$.

Notice that, for a fixed admissible control $u_{\ell}$ from the measurable set\, $\mathcal{U}_{\ell}$, if we obtain a representation for Equation~\eqref{Eq5} as a minimal cost for an associated stochastic optimal control problem, then we will be able to obtain a representation for $V_{\ell, \theta_{\ell}}^{\epsilon}$ as a value function for a stochastic differential game. This further allow us to link this progressively measurable control $u_{\ell}$ in the original control problem with a strategy for the maximizing player of the associated stochastic differential game. Furthermore, such a connection between the risk-sensitive value function and a deterministic differential game can be made immediately, when the small random perturbation vanishes in the limit. 

Before concluding this section, it is worth mentioning that some interesting studies on risk-sensitive control problem for dynamical systems with small random perturbations have been reported in literature (for example, see \cite{DupMc97} using PDE viscosity solution techniques; see \cite{BouDu01} using the probabilistic argumentation and the variational representation for degenerate diffusion processes; see also \cite{DaiMR96}, \cite{ElHam03} or \cite{FleMc95} for some connections between the risk-sensitive stochastic control and dynamic games).
 
An outline of the paper is as follows. In Section~\ref{S2}, we introduce a family of two-player differential games -- where {\it Player}-$1$ will attempt to maximize the mean escape time criterion corresponding to each of the subsystems; while {\it Player}-$2$ will attempt to minimize it. In this section, we also provide some preliminary results that are useful for proving our main results. In Section~\ref{S3}, we present our main results -- where we consider a risk-sensitive version of the mean escape time criterion with respect to each of the subsystems. Using the variational representation, we characterize the risk-sensitive escape control for the diffusion process as the lower and upper values of the associated stochastic differential game. Finally, we comment on the implication of our results, where one is also interested in evaluating the performance of the risk-sensitive escape control for the diffusion process, when there is some norm-bounded modeling error in the distributed control system.

\section{Preliminary Results} \label{S2}

\subsection{A Differential Game Formalism} \label{S2(1)}

In this subsection, we consider a family of two-player differential games. For a fixed $\ell \in \{2, 3, \ldots, n\}$, at each time $t \in [0,\, \infty)$, {\it Player}-$1$ picks a strategy $u_{\ell}(t)$ from the admissible control space $\mathcal{U}_{\ell}$, and {\it Player}-$2$ picks a control $v_{\ell}(t)$ from $\mathbb{R}^m$ in such a way that the functions $t \mapsto u_{\ell}(t)$ and $t \mapsto v_{\ell}(t)$ belong to the strategy sets
\begin{align}
M_{\ell} = \biggl\{u_{\ell} \colon [0,\, \infty) \rightarrow \mathcal{U}_{\ell} \, \bigl \vert \,u_{\ell} \,\,\text{is progressively measurable} \biggr\}   \label{Eq8}
\end{align}
and
\begin{align}
N_{\ell} = \biggl\{v_{\ell} \colon [0,\, \infty) \rightarrow \mathbb{R}^m \, \biggl \vert \,  \int_0^T  \bigl\vert v_{\ell}(t) \bigr \vert^2 dt  < \infty \quad \forall T < \infty \biggr\}, \label{Eq9}
\end{align}
respectively. Here, we also identify that $M_{\ell}$ and $N_{\ell}$, for any $\ell \in \{2, 3, \ldots, n \}$, as metric spaces under any metric which is equivalent to convergence in $\mathscr{L}^2\bigl([0,\,\ T], \mathbb{R}^{r_i} \bigr)$ and  $\mathscr{L}^2\bigl([0,\,\ T], \mathbb{R}^m \bigr)$.

Suppose that both players have played the game up to the $(\ell-1)$th-stage (see Footnote~\ref{FN}). Let $u_{i}^{\ast}(t) \in \mathcal{U}_i$, $\forall t \in [0,\, \infty)$, for $i =1, 2, \ldots, \ell-1$, be the admissible control strategies picked by the maximizing {\it Player}-$1$. Then, at the $\ell$th-stage, the dynamics of the game is given by the following differential equations
\begin{align}
\left.\begin{array}{l}
d x^{1}(t) = m_1\bigl(x^{1}(t), u_1^{\ast}(t)\bigr) dt + \sigma\bigl(x^{1}(t)\bigr)v_{\ell}(t) dt \\
d x^{2}(t) = m_2\bigl(x^{1}(t), x^{2}(t), u_2^{\ast}(t)\bigr) dt  \\
 \quad\quad\quad ~ \vdots  \\
d x^{\ell-1}(t) = m_{\ell-1}\bigl(x^{1}(t), x^{2}(t), \ldots, x^{\ell-1}(t), u_{\ell-1}^{\ast}(t)\bigr) dt  \\
d x^{\ell}(t) = m_{\ell}\bigl(x^{1}(t), x^{2}(t), \ldots, x^{\ell}(t), u_{\ell}(t)\bigr) dt\\
\quad \quad \quad \quad \quad \quad\quad \quad \quad  \bigl(x_0^{1}, x_0^{2}, \ldots, x_0^{\ell}\bigr) \in \Omega_{\ell}, \,\, t \ge 0
 \end{array}\right\}  \label{Eq10}
\end{align}
with an associated cost criterion
\begin{align}
J_{\ell} \bigl(u_{\ell}, v_{\ell}\bigr) = \frac{1}{2} \int_0^{\tau_{\ell}^0} \bigl\vert v_{\ell}(t) \bigr \vert^2 dt + \theta_{\ell} \tau_{\ell}^0, \label{Eq11}
\end{align}
where $\tau_{\ell}^0 \triangleq \inf \bigl\{ t > 0 \, \vert \, x^{\ell}(t) \notin D_{\ell}  \bigr\}$. Note that the goal of {\it Player}-$1$ is to maximize $J_{\ell}$ with respect to $u_{\ell}$ and while that of {\it Player}-$2$ is to minimize it with respect to $v_{\ell}$, for each $\ell \in \{2, 3, \ldots, n \}$. Here, we remark that {\it Player}-$1$ will attempt preventing the diffusion process $x^{\ell}(t)$ from leaving the given domain $D_{\ell}$ (i.e., representing the exact control in risk-sensitive problem); while {\it Player}-$2$ will attempt forcing out the diffusion process from the domain (i.e., acting the role of the disturbance in the distributed control system).\footnote{We remark that the admissible control strategy $u_{\ell}$ picked by {\it Player}-$1$ affects only the dynamics of the game, not directly the cost criterion.}

Furthermore, a mapping $\alpha_{\ell} \colon N_{\ell} \rightarrow M_{\ell}$ is said to be a strategy for the maximizing player if it is measurable and, for $v_{\ell},\, \hat{v}_{\ell} \in N_{\ell}$,
\begin{align*}
v_{\ell}(t) = \hat{v}_{\ell}(t), \quad \forall t \in [0,\, s]
\end{align*}
implies
\begin{align*}
\alpha_{\ell}[v_{\ell}](t) = \alpha_{\ell}[\hat{v}_{\ell}](t), \quad \forall t \in [0,\, s],
\end{align*}
almost everywhere, for every $s \in [0,\, \infty)$. 

Similarly, a mapping $\beta_{\ell} \colon M_{\ell} \rightarrow N_{\ell}$ is a strategy for the minimizing player if it is measurable and, for $u_{\ell},\, \hat{u}_{\ell} \in M_{\ell}$,
\begin{align*}
u_{\ell}(t) = \hat{u}_{\ell}(t), \quad \forall t \in [0,\, s]
\end{align*}
implies
\begin{align*}
\beta_{\ell}[u_{\ell}](t) = \beta_{\ell}[\hat{u}_{\ell}](t), \quad \forall t \in [0,\, s],
\end{align*}
almost everywhere, for every $s \in [0,\, \infty)$.\footnote{\label{FN}Note that during each expanding construction (i.e., when a new subsystem is  added to the existing distributed control system), we assume that both players play a differential game. For example, for $\ell = 2$, the dynamics of the game is given by
\begin{align*}
\left.\begin{array}{l}
d x^{1}(t) = m_1\bigl(x^{1}(t), u_1^{\ast}(t)\bigr) dt + \sigma\bigl(x^{1}(t)\bigr)v_2(t) dt \\
d x^{2}(t) = m_2\bigl(x^{1}(t), x^{2}(t), u_2(t)\bigr) dt  \\
 \quad \quad\quad \quad \quad  \bigl(x_0^{1}, x_0^{2}\bigr) \in \Omega_{2} = D_1 \times D_2, \,\, t \ge 0,
 \end{array}\right\}
\end{align*}
with an associated cost criterion
\begin{align*}
J_{2} \bigl(u_{2}, v_{2}\bigr) = \frac{1}{2} \int_0^{\tau_{2}^0} \bigl\vert v_{2}(t) \bigr \vert^2 dt + \theta_{2} \tau_{2}^0
\end{align*}
and an exit-time $\tau_{2}^0 \triangleq \inf \bigl\{ t > 0 \, \vert \, x^{2}(t) \notin D_{2}  \bigr\}$ such that {\it Player}-$1$ optimally picks a strategy (in the sense of best-response correspondence) to the stagey of {\it Player}-$2$. Then, the game advances to the next stage, i.e., $\ell = 3$, and continues until $n$.}

Let us denote the set of all maximizing strategies by $\Gamma_{\ell}$ and the set of all minimizing strategies by $\Lambda_{\ell}$. Furthermore, let us define the lower and the upper values of the differential game at the $\ell$th-stage by
\begin{align}
I_{\ell}^{-}\bigl(x_0^{1}, x_0^{2}, \ldots, x_0^{\ell}\bigr) = \inf_{\beta_{\ell} \in \Lambda_{\ell}} \sup_{u_{\ell} \in \Gamma_{\ell}} J_{\ell} \bigl(u_{\ell}, v_{\ell}\bigr) \label{Eq12}
\end{align}
and
\begin{align}
I_{\ell}^{+}\bigl(x_0^{1}, x_0^{2}, \ldots, x_0^{\ell}\bigr) = \sup_{u_{\ell} \in \Gamma_{\ell}}  \inf_{\beta_{\ell} \in \Lambda_{\ell}} J_{\ell} \bigl(u_{\ell}, v_{\ell}\bigr),  \label{Eq13}
\end{align}
for each ${\ell} \in \{2, 3, \ldots, n\}$, respectively. Moreover, if 
\begin{align}
I_{\ell}^{-}\bigl(x_0^{1}, x_0^{2}, \ldots, x_0^{\ell}\bigr) = I_{\ell}^{+}\bigl(x_0^{1}, x_0^{2}, \ldots, x_0^{\ell}\bigr),  \label{Eq14}
\end{align}
then the differential game has a value.

\begin{remark} \label{R3}
Note that the greatest payoff that {\it Player}-$1$ (i.e., the maximizing player) can force is called a lower value of the game and, similarly, the least value that {\it Player}-$2$ (i.e., the minimizing player) can force is termed an upper value of the game. In Section~\ref{S3}, we provide conditions under which these values coincide.
\end{remark}

\subsection{Additional Preliminary Results} \label{S2(3)}
In this subsection, we provide additional results that will be useful for proving our main results in Section~\ref{S3}.

\begin{definition} \label{D1}
We define $\mathscr{V}_i$ to be the set of all $\mathbb{R}^m$-valued $\mathscr{F}_t$-progressively measurable processes $v_i \triangleq \bigl\{v_i(t),\, t \in[0,\,\infty)  \bigr \}$, that satisfies
\begin{align*}
\mathbb{E}_{x_{\widehat{1,\ell}}}^{\epsilon} \biggl\{ \int_0^T  \bigl\vert v_i(t) \bigr \vert^2 dt \biggr \} < \infty \quad \forall T < \infty,
\end{align*}
for each $i = 1, 2, \ldots, n$.
\end{definition}

\begin{lemma} \label{L1} [Variational representation formula (cf. \cite[Proposition~2.5 or Theorem~5.1]{BouDu98})]
For a fixed $T < \infty$, let $f \colon \mathscr{C} \bigl([0,\,T], \mathbb{R}^m \bigr) \rightarrow \mathbb{R}$ be any Borel measurable bounded function. Then
\begin{align}
-\log \mathbb{E}_{x_{\widehat{1,\ell}}}^{\epsilon} \exp \biggl\{-f\bigl(W\bigr)\biggr\}  = \inf_{v_{\ell} \in \mathscr{V}_{\ell}} \mathbb{E}_{x_{\widehat{1,\ell}}}^{\epsilon} \biggl\{ \frac{1}{2} \int_0^T \bigl\vert v_{{\ell}}(s) \bigr \vert^2ds + f\biggl(W + \int_0^{\cdot} v_{{\ell}}(s) ds \biggr)\biggr\}, \label{Eq15}
\end{align}
for any ${\ell} \in \{2, 3, \ldots, n\}$. 
\end{lemma}

\begin{lemma} \label{L2}
For any $v_{\ell} \in \mathscr{V}_{\ell}$, with ${\ell} \in \{2, 3, \ldots, n\}$, and $T < \infty$, let $\mu^{v_{\ell}}$ be a measure induced on $\mathscr{C} \bigl([0,\,T], \mathbb{R}^m\bigr)$ by $W + \int_0^{\cdot} v_{\ell}(t) dt$ under $P$. Then, the relative entropy of $\mu^{v_{\ell}}$ with respect to $P$ satisfies the following 
\begin{align}
\mathcal{R}\bigl(\mu^{v_{\ell}} \bigl \Vert P\bigr) \le \mathbb{E}_{x_{\widehat{1,\ell}}}^{\epsilon} \biggl\{ \int_0^T  \bigl\vert v_{\ell}(t) \bigr \vert^2 dt \biggr \}.  \label{Eq16}
\end{align}
\end{lemma}

Let $\mathscr{S}$ a Polish space (i.e., a complete separable metric space), with a Borel $\sigma$-algebra, and let $\mathscr{P}(\mathscr{S})$ be the set of measures defined on $\mathscr{S}$ that satisfies the usual hypotheses (e.g., see \cite{Kry80}).
\begin{lemma} \label{L3} [cf. \cite[Theorem~2.1]{Bill68}]
Consider a sequence of measures $\bigl\{\mu_{k},\, k \in \mathbb{N} \bigr\}$  in $\mathscr{P}(\mathscr{S})$ satisfying
\begin{align}
\sup_{k \in \mathbb{N}} \mathcal{R}\bigl(\mu_k \bigl \Vert P\bigr) < \infty, \label{Eq17}
\end{align}
where $P \in \mathscr{P}(\mathscr{S})$. Let $f \colon \mathscr{S} \rightarrow \mathbb{R}$ be a Borel-measurable function. Then, the followings hold
\begin{enumerate} [(i)]
\item  if $\mu_k$ weakly converges to another measure $\mu$ as $k \rightarrow \infty$, then
\begin{align}
\lim_{k \rightarrow \infty} \int_{\mathscr{S}} f d \mu_k =  \int_{\mathscr{S}} f d \mu, \label{Eq18-a}
\end{align}
\item if $\bigl\{ f_k,\, k \in \mathbb {N} \bigr\}$ is a sequence of uniformly bounded functions that almost surely converges to $f$, then
\begin{align}
\lim_{k \rightarrow \infty} \int_{\mathscr{S}} f_k d \mu_k =  \int_{\mathscr{S}} f d \mu. \label{Eq18-b}
\end{align}
 \end{enumerate}
\end{lemma}

\section{Main Results} \label{S3}
\subsection{Risk-Sensitive Escape Control Problem} \label{S3(1)}
In this subsection, we relate the lower and upper values of the associated differential game with the risk-sensitive escape control problem for the diffusion process. In particular, using the variational representation (e.g, see \cite{BouDu98} or \cite{DupMc97}), we present our main results, i.e., Proposition~\ref{P1} and Proposition~\ref{P2}.

For each fixed admissible control $u_{\ell} \in \mathcal{U}_{\ell}$, the following proposition (which is a direct consequence of Lemma~\ref{L1}) characterizes the risk-sensitive escape control problem (cf. Equations~\eqref{Eq5} and \eqref{Eq7}) with an associated stochastic differential game (cf. Equations~\eqref{Eq22} below).

\begin{proposition} \label{P1}
Suppose that, for a fixed ${\ell} \in \{2, 3, \ldots, n\}$, the admissible optimal controls $u_{i}^{\ast}(t) \in \mathcal{U}_i$, $\forall t \in [0,\, \infty)$, for $i =1, 2, \ldots, \ell-1$, are given. Consider any admissible control $u_{\ell} \in \mathcal{U}_{\ell}$. Further, for every $s \in [0,\, \infty)$, let $\gamma_{\ell}^s$ be a Borel measurable function such that $u_{\ell}(t)=\gamma_{\ell}^s\bigl(t, W(r), 0\le r \le t \bigr)$ for $t \in [0,\, s]$, w.p.1. Let the exit-time $\tau_{\ell}^{\epsilon}$ be given by 
\begin{align*}
\tau_{\ell}^{\epsilon} = \inf \bigl\{ t > 0 \, \bigl\vert \, x^{\ell}(t) \notin D_{\ell} \bigr\},
\end{align*}
which is associated with the following diffusion processes $\bigl(x^1, x^2, \ldots, x^{\ell}\bigr)$, i.e.,
\begin{align}
\left.\begin{array}{l}
d x^{1}(t) = m_1\bigl(x^{1}(t), u_1^{\ast}(t)\bigr) dt + \sqrt{\epsilon} \sigma\bigl(x^{1}(t)\bigr)dW(t) \\
d x^{2}(t) = m_2\bigl(x^{1}(t), x^{2}(t), u_2^{\ast}(t)\bigr) dt  \\
 \quad\quad\quad ~ \vdots  \\
d x^{\ell-1}(t) = m_{\ell-1}\bigl(x^{1}(t), x^{2}(t), \ldots, x^{\ell-1}(t), u_{\ell-1}^{\ast}(t)\bigr) dt  \\
d x^{\ell}(t) = m_{\ell}\bigl(x^{1}(t), x^{2}(t), \ldots, x^{\ell}(t), u_{\ell}(t)\bigr) dt\\
  \quad \quad \quad \quad \quad \quad\quad \quad \quad  \bigl(x_0^{1}, x_0^{2}, \ldots, x_0^{\ell}\bigr) \in \Omega_{\ell}, \,\, t \ge 0
\end{array}\right\}  \label{Eq19}
\end{align}
Then, the following variational representation holds
\begin{align}
-\epsilon \log \mathbb{E}_{x_{\widehat{1,\ell}}}^{\epsilon} \exp \biggl\{-\frac{1}{\epsilon} \theta_{\ell}\, \tau_{\ell}^{\epsilon} \biggr\} = \inf_{v_{\ell} \in \mathscr{V}_{\ell}} \mathbb{E}_{x_{\widehat{1,\ell}}}^{\epsilon} \biggl\{ \frac{1}{2} \int_0^{\tilde{\tau}_{\ell}^{\epsilon}} \bigl\vert v_{\ell}(s) \bigr \vert^2ds + \theta_{\ell} \tilde{\tau}_{\ell}^{\epsilon} \biggr\}, \label{Eq20}
\end{align}
where the exit-time $\tilde{\tau}_{\ell}^{\epsilon}$ is given by
\begin{align}
\tilde{\tau}_{\ell}^{\epsilon} = \inf \Bigl\{ t > 0 \, \bigl\vert \, \tilde{x}^{\ell}(t) \notin D_{\ell} \Bigr\}, \label{Eq21}
\end{align}
which is associated with the following diffusion processes $\bigl(\tilde{x}^1, \tilde{x}^2, \ldots, \tilde{x}^{\ell}\bigr)$, i.e.,
\begin{align}
\left.\begin{array}{l}
d \tilde{x}^{1}(t) = m_1\bigl(\tilde{x}^{1}(t), \tilde{u}_1^{\ast}(t)\bigr) dt + \sigma\bigl(\tilde{x}^{1}(t))v_{\ell}(t) dt + \sqrt{\epsilon} \sigma\bigl(\tilde{x}^{1}(t))dW(t) \\
d \tilde{x}^{2}(t) = m_2\bigl(\tilde{x}^{1}(t), \tilde{x}^{2}(t), \tilde{u}_2^{\ast}(t)\bigr) dt  \\
 \quad\quad\quad ~ \vdots  \\
d \tilde{x}^{\ell-1}(t) = m_{\ell-1}\bigl(\tilde{x}^{1}(t), \tilde{x}^{2}(t), \ldots, \tilde{x}^{\ell-1}(t), \tilde{u}_{\ell-1}^{\ast}(t)\bigr) dt  \\
d \tilde{x}^{\ell}(t) = m_{\ell}\bigl(\tilde{x}^{1}(t), \tilde{x}^{2}(t), \ldots, \tilde{x}^{\ell}(t), \tilde{u}_{\ell}(t)\bigr) dt\\
 \bigl(\tilde{x}^{1}(0), \tilde{x}^{2}(0), \ldots, \tilde{x}^{\ell}(0)\bigr) \triangleq \bigl(x_0^{1}, x_0^{2}, \ldots, x_0^{\ell}\bigr) \in \Omega_{\ell}, \,\, t \ge 0
\end{array}\right\}  \label{Eq22}
\end{align}
Moreover, the admissible control $\tilde{u}_{\ell}$ satisfies
\begin{align}
\tilde{u}_{\ell}(t) = \gamma_{\ell}^s\biggl(t, W(r) + \frac{1}{\sqrt{\epsilon}}\int_0^{r} v_{\ell}(z) dz, 0 \le r \le t \biggr), \quad \forall t \in [0,\,s], \quad \text{w.p.1}, \label{Eq23}
\end{align}
for any $s \in [0,\,\infty)$.
\end{proposition}
The following proposition provides conditions under which the lower and upper values of the associated differential game (i.e., quantities in Equations~\eqref{Eq12} and \eqref{Eq13}) will coincide.
\begin{proposition} \label{P2}
For a fixed ${\ell} \in \{2, 3, \ldots, n\}$, let $I_{\ell}^{-}\bigl(x_0^{1}, x_0^{2}, \ldots, x_0^{\ell}\bigr)$ and $I_{\ell}^{+}\bigl(x_0^{1}, x_0^{2}, \ldots, x_0^{\ell}\bigr)$ be the lower and upper values of the associated differential game given in Equations~\eqref{Eq12} and \eqref{Eq13}. Suppose that $u_{i}^{\ast}(t) \in \mathcal{U}_i$, $\forall t \in [0,\, \infty)$, for $i =1, 2, \ldots, \ell-1$, are admissible optimal controls. For any $u_{\ell} \in \mathcal{U}_{\ell}$, let $\bigl(x^1, x^2, \ldots, x^{\ell}\bigr)$ be the unique solution to Equation~\eqref{Eq19}. Then, 
\begin{enumerate} [(i)]
\item the lower value of the game satisfies
\begin{align}
\limsup_{\epsilon \rightarrow 0} \sup_{u_{\ell} \in \mathcal{U}_{\ell}} -\epsilon \log \mathbb{E}_{x_{\widehat{1,\ell}}}^{\epsilon} \exp \biggl\{-\frac{1}{\epsilon} \theta_{\ell}\, \tau_{\ell}^{\epsilon} \biggr\} \le I_{\ell}^{-}\bigl(x_0^{1}, x_0^{2}, \ldots, x_0^{\ell}\bigr), \label{Eq24}
\end{align}
\item for a given $\kappa > 0$, there exists a measurable function $g_{\ell} \colon \mathscr{C} \bigl([0,\,T], \mathbb{R}^m\bigr) \rightarrow M_{\ell}$ such that the upper value of the game satisfies
\begin{align}
\liminf_{\epsilon \rightarrow 0} \sup_{u_{\ell} \in \mathcal{U}_{\ell}} -\epsilon \log \mathbb{E}_{x_{\widehat{1,\ell}}}^{\epsilon} \exp \biggl\{-\frac{1}{\epsilon} \theta_{\ell}\, \tau_{\ell}^{\epsilon} \biggr\} \ge I_{\ell}^{+}\bigl(x_0^{1}, x_0^{2}, \ldots, x_0^{\ell}\bigr) - \kappa, \label{Eq25}
\end{align}
with $u_{\ell} = g_{\ell}[\sqrt{\epsilon} W]$ (i.e., when the maximizing player picks such a strategy),
\item if the lower and upper values of the game coincides, i.e.,
\begin{align}
I_{\ell}^{-}\bigl(x_0^{1}, x_0^{2}, \ldots, x_0^{\ell}\bigr) = I_{\ell}^{+}\bigl(x_0^{1}, x_0^{2}, \ldots, x_0^{\ell}\bigr), \label{Eq26}
\end{align}
then the game has a value.
\end{enumerate} 
\end{proposition}

\subsection{Remarks on the Robust Analysis Problem} \label{S3(2)}
In this subsection, we briefly remark on the implication of our main results -- where one is also interested in evaluating the robust performance of the risk-sensitive escape control, when there is some norm-bounded modeling error in the distributed control system (see \cite{DupJaPe00} for related discussion, but in different context).

In what follows, we assume that that the statements in Propositions~\ref{P1} and \ref{P2} are true. Suppose that, for a fixed ${\ell} \in \{2, 3, \ldots, n\}$, $u_{i}^{\ast}(t) \in \mathcal{U}_i$, $\forall t \in [0,\, \infty)$, for $i =1, 2, \ldots, \ell$, are the admissible optimal control strategies picked by {\it Player}-$1$. Further, we consider the following distributed control system (which contains $\ell$ subsystems)
\begin{align}
\left.\begin{array}{l}
d \hat{x}^{1}(t) = m_1\bigl(\hat{x}^{1}(t), {u}_{1}^{\ast}(t)\bigr) dt + \sigma\bigl(\hat{x}^{1}(t)) \hat{v}_{\ell}(t) dt \\
d \hat{x}^{2}(t) = m_2\bigl(\hat{x}^{1}(t), {u}_{2}^{\ast}(t), \hat{u}_{2}(t) \bigr) dt  \\
 \quad\quad\quad ~ \vdots  \\
d \hat{x}^{\ell}(t) = m_{\ell}\bigl(\hat{x}^{1}(t), \hat{x}^{2}(t), \ldots, \hat{x}^{\ell}(t), {u}_{\ell}^{\ast}(t)\bigr) dt\\
\quad  \bigl(\hat{x}^{1}(0), \hat{x}^{2}(0), \ldots, \hat{x}^{\ell}(0)\bigr) \triangleq \bigl(x_0^{1}, x_0^{2}, \ldots, x_0^{\ell}\bigr) \in \Omega_{\ell}, \,\, t \ge 0
\end{array}\right\} \label{Eq27}
\end{align}
where $\hat{v}_{\ell}(t) \in N_{\ell}$.

Define the value function $V_{\ell, \theta_{\ell}}^0 = V_{\ell, \theta_{\ell}}^0\bigl(x_0^{1}, x_0^{2}, \ldots, x_0^{\ell}\bigr)$ as
\begin{align}
V_{\ell, \theta_{\ell}}^0 \triangleq \inf_{\hat{v}_{\ell} \in N_{\ell}} \biggl \{ \frac{1}{2} \int_0^{\tau_{\ell}^0} \bigl\vert \hat{v}_{\ell}(t) \bigr \vert^2 dt + \theta_{\ell} \tau_{\ell}^0 \biggr\}, \label{Eq28}
\end{align}
where $\tau_{\ell}^0 = \inf \bigl\{ t > 0 \, \vert \, x^{\ell}(t) \notin D_{\ell}  \bigr\}$, with $\theta_{\ell} > 0$. Then, for any $\hat{v}_{\ell} \in N_{\ell}$, with $\bigl \Vert \hat{v}_{\ell} \bigr \Vert_{\infty} < \infty$, we have the following inequalities
\begin{align}
V_{\ell, \theta_{\ell}}^0 &\le \frac{1}{2} \int_0^{\tau_{\ell}^0} \bigl\vert \hat{v}_{\ell}(t) \bigr \vert^2 dt + \theta_{\ell} \tau_{\ell}^0, \notag \\ 
                                 & \le \frac{\tau_{\ell}^0}{2} \bigl \Vert \hat{v}_{\ell} \bigr \Vert_{\infty}^2 + \theta_{\ell} \tau_{\ell}^0 \label{Eq29}
\end{align}
and
\begin{align}
\tau_{\ell}^0 \ge \frac{V_{\ell, \theta_{\ell}}^0}{\frac{1}{2} \bigl \Vert \hat{v}_{\ell} \bigr \Vert_{\infty}^2 + \theta_{\ell}}.  \label{Eq30}
\end{align}
Suppose that $\hat{v}_{\ell}(t) \triangleq \Delta m_{\ell}(t)$, where $\Delta m_{\ell}(t)$ is interpreted as a modeling error in Equation~\eqref{Eq26}. Further, assume that the value $\tau_{\ell}^0$ is used as a qualitative measure on the performance of the distributed control system. For a given specification $L_{\ell}$, with $\tau_{\ell}^0 \ge L_{\ell}$, if there exists a design parameter $\theta_{\ell}^{\ast} > 0$ such that
\begin{align}
\frac{V_{\ell, \theta_{\ell}^{\ast}}^0}{\theta_{\ell}^{\ast}} > L_{\ell}. \label{Eq31}
\end{align}
Then, we obtain an upper bound on the norm of the modeling error, which guarantees the desired performance against all modeling errors satisfying such a norm bound. That is, if 
\begin{align}
\bigl \Vert \Delta m_{\ell}\bigr \Vert_{\infty}^2 \le 2\left[ \frac{V_{\ell, \theta_{\ell}^{\ast}}^0}{L_{\ell}} - \theta_{\ell}^{\ast} \right], \label{Eq32}
\end{align}
then we have
\begin{align}
\tau_{\ell}^0 \ge \frac{V_{\ell, \theta_{\ell}^{\ast}}^0}{\frac{1}{2} \bigl \Vert \Delta m_{\ell} \bigr \Vert_{\infty}^2 + \theta_{\ell}^{\ast}} \ge L_{\ell}. \label{Eq33}
\end{align}
Moreover, the above equation (together with Equation~\eqref{Eq6}) further implies the following
\begin{align}
 \frac{V_{i, \theta_{i}^{\ast}}^0}{\frac{1}{2} \bigl \Vert \Delta m_{i} \bigr \Vert_{\infty}^2 + \theta_{i}^{\ast}} \ge \frac{V_{i+1, \theta_{i+1}^{\ast}}^0}{\frac{1}{2} \bigl \Vert \Delta m_{i+1} \bigr \Vert_{\infty}^2 + \theta_{i+1}^{\ast}}, \label{Eq34}
\end{align}
for $i =1, 2, \ldots, \ell-1$, with $\ell \in \{2, 3, \ldots, n\}$.
\begin{remark} \label{R4}
Note that, for each $i = 1, 2, \ldots, n$, the norm on the modeling error is inversely proportional to the design specification $L_{i}$, and, therefore, the robustness of the distributed control system increases as the bound on the performance measure decreases.
\end{remark}

\end{document}